\documentclass[reqno,11pt]{amsart} 
\usepackage{amsmath}
\usepackage{amssymb}
\usepackage{amsthm}

\newcommand\NoBlackBoxes{\global\overfullrule0pt}
\NoBlackBoxes
\theoremstyle{plain}

\begin{document}

\title{H\"offding's kernels and periodic \\
covariance representations}

\author{Sergey G. Bobkov$^{1}$}
\thanks{1) 
School of Mathematics, University of Minnesota, Minneapolis, MN, USA
}

\author{Devraj Duggal$^{2}$}
\thanks{2) School of Mathematics, University of Minnesota, Minneapolis, MN, USA
}


\begin{abstract}
We start with a brief survey on H\"offding's kernels, its properties, related
spectral decompositions, and discuss marginal distributions of H\"offding measures.
In the second part of this note, one-dimensional covariance representations 
are considered over compactly supported probability distributions in the class 
of periodic smooth functions. H\"offding's kernels are used in the construction
of mixing measures whose marginals are multiples of given probability distributions,
leading to optimal kernels in periodic covariance representations.
\end{abstract} 
 
\maketitle
\markboth{Sergey G. Bobkov and Devraj Duggal}{Periodic covariance representations}

\def\theequation{\thesection.\arabic{equation}}
\def\E{{\mathbb E}}
\def\R{{\mathbb R}}
\def\C{{\mathbb C}}
\def\P{{\mathbb P}}
\def\Z{{\mathbb Z}}
\def\S{{\mathbb S}}
\def\I{{\mathbb I}}
\def\T{{\mathbb T}}

\def\s{{\mathbb s}}

\def\G{\Gamma}

\def\Ent{{\rm Ent}}
\def\var{{\rm Var}}
\def\Var{{\rm Var}}
\def\cov{{\rm cov}}
\def\V{{\rm V}}

\def\H{{\rm H}}
\def\Im{{\rm Im}}
\def\Tr{{\rm Tr}}
\def\s{{\mathfrak s}}
\def\A{{\mathfrak A}}
\def\m{{\mathfrak m}}

\def\k{{\kappa}}
\def\M{{\cal M}}
\def\Var{{\rm Var}}
\def\Ent{{\rm Ent}}
\def\O{{\rm Osc}_\mu}

\def\ep{\varepsilon}
\def\phi{\varphi}
\def\vp{\varphi}
\def\F{{\cal F}}

\def\be{\begin{equation}}
\def\en{\end{equation}}
\def\bee{\begin{eqnarray*}}
\def\ene{\end{eqnarray*}}

\thispagestyle{empty}

\section{{\bf Generalized H\"offding's formula}}
\setcounter{equation}{0}

\vskip2mm
\noindent
Given two random variables $X$ and $Y$, 
the generalized H\"offding's covariance formula indicates that, for all 
``regular" functions $u$ and $v$ on the real line,
\be
\cov(u(X),v(Y)) = \int_{-\infty}^\infty \int_{-\infty}^\infty
u'(x) v'(y)\,H(x,y)\,dx\,dy,
\en
where
\be
H(x,y) = \P\{X \leq x, Y \leq y\} - \P\{X \leq x\}\, 
\P\{Y \leq y\}, \quad x,y \in \R.
\en
The case of the identical functions $f(x)=x$ and $g(y)=y$ corresponds to
H\"offding \cite{H} (provided that $X$ and $Y$ finite second moments). 
The history of this remarkable identity may be found
in Lo \cite{Lo}, together with generalizations and refinements of the
previous results by Mardia \cite{M}, Sen \cite{Sen}, Cuadras \cite{C1,C2};
see also recent works by Saumard and Wellner \cite{S-W1, S-W2}.

Block and Fang \cite{B-F} proposed an extension of the original H\"offding's 
formula to more than two variables. Let us however restrict ourselves to 
the particular case of (1.1) with $X = Y$ and rewrite this relation 
as a covariance identity with respect to the distribution $\mu$ of $X$:
\be
{\rm cov}_\mu(u,v) = \int_{-\infty}^\infty \int_{-\infty}^\infty 
u'(x) v'(y)\,d\lambda(x,y).
\en
Here, one may require that $\lambda$ be a positive, locally finite measure on 
the plane $\R \times \R$ (that is, finite on compact sets). According to (1.1)-(1.2), 
the identity (1.3) holds true, when $\lambda$ 
is absolutely continuous over the Lebesgue measure with density
\be
H_\mu(x,y) = \frac{d\lambda(x,y)}{dx\,dy} = F(x \wedge y)\,(1 - F(x \vee y)), \quad
x,y \in \R,
\en
where $F(x) = \P\{X \leq x\} = \mu((-\infty,x])$
is the associated distribution function. We adopt the standard notations
$x \wedge y = \min(x,y)$, $x \vee y = \max(x,y)$. 

\vskip5mm
{\bf Definition.} We call $\lambda = \lambda_\mu$ the H\"offding measure
and its density $H = H_\mu$ the H\"offding kernel associated to $\mu$.

\vskip2mm
For example, if $\mu = p\delta_a + q\delta_b$ is the Bernoulli measure 
assigning the weights $p \in (0,1)$ and $q=1-p$ to the points $a < b$, then
$\lambda_\mu = pq U$ where $U$ is the uniform distribution on the square
$(a,b) \times (a,b)$.

Let us state this consequence of (1.1)-(1.2) once more in the next statement 
with emphasis on the uniqueness part in the representation (1.3). 
As the weakest requirement, one may consider the latter identity 
in the class $C^\infty_b$ of all
functions $u,v:\R \rightarrow \R$ having $C^\infty$-smooth, compactly
supported derivatives (in which case $u$ and $v$ are bounded).

\vskip5mm
{\bf Theorem 1.1.} {\sl Given a probability measure $\mu$ on the real line, 
$(1.3)$ holds true for all $u,v \in C^\infty_b$ with a unique positive, 
locally finite measure
$\lambda = \lambda_\mu$. Moreover, $(1.3)$ extends to all locally 
absolutely continuous, complex-valued functions $u,v$ such that
\be
\int_{-\infty}^\infty \int_{-\infty}^\infty 
|u'(x)|\, |u'(y)|\,d\lambda(x,y) < \infty, \quad
\int_{-\infty}^\infty \int_{-\infty}^\infty 
|v'(x)|\, |v'(y)|\,d\lambda(x,y) < \infty,
\en
where the derivatives $u'$ and $v'$ are understood in the Radon-Nikodym sense. 
The H\"offding measure $\lambda$ is finite, if and only if
$\mu$ has a finite second moment.
}

\vskip5mm
The condition (1.5) insures that the function $u'(x) v'(y)$ is integrable 
over $\lambda$, and also implies that $u(X)$ and $v(X)$ have finite second 
moments. Hence both sides in (1.3) are well-defined and finite. For the sake 
of completeness, we sketch a short proof of Theorem 1.1 and give a few 
remarks on the existing (slightly different) formulations in the end (Section 10). 

Let us note that, being applied with $u(x) = v(x) = x$, the identity (1.3) 
shows that the total mass of the H\"offding measure is the variance
$$
\lambda(\R \times \R) = \int_{-\infty}^\infty \int_{-\infty}^\infty 
H(x,y)\,dx\,dy = \Var(X).
$$
Once this measure is finite, it may also be described via its
Fourier-Stieltjes transform in terms of the characteristic function
of the random variable $X$,
$$
f(t) = \E\,e^{itX} = \int_{-\infty}^\infty e^{itx}\,d\mu(x),
\quad t \in \R.
$$
Namely, applying (1.3) to the exponential functions $u(x) = e^{itx}$ 
and $v(y) = e^{isy}$ ($t,s \in \R$), we obtain an explicit formula
\be
\hat \lambda(t,s) = \int_{-\infty}^\infty \int_{-\infty}^\infty 
e^{itx + isy}\,d\lambda(x,y) =
\frac{f(t) f(s) - f(t+s)}{ts} \quad (t,s \neq 0).
\en
In particular, this provides the uniqueness part in Theorem 1.1.

Thus, the expression on the right-hand side of (1.6) represents a positive 
definite function in two real variables, as long as the characteristic 
function $f$ is twice differentiable. Moreover, 
there is a similar property of H\"offding's kernels themselves.

\vskip5mm
\section{{\bf Positive definiteness of H\"offding's kernels}}
\setcounter{equation}{0}

\vskip2mm
\noindent
Let $X$ be a random variable with a non-degenerate distribution $\mu$,
so that the associated H\"offding measure $\lambda$ is non-zero.
From (1.5) we have a variance representation
\be
\Var(u(X)) = \int_{-\infty}^\infty \int_{-\infty}^\infty 
u'(x) u'(y)\,d\lambda(x,y),
\en
and the substitution $f = u'$ leads to
\be
\int_{-\infty}^\infty \int_{-\infty}^\infty f(x) f(y) H(x,y)\,dx\,dy
 \geq 0,
\en
which holds for any measurable function $f$ on the real line such that
the integral is well-defined in the Lebesgue sense. In other words:

\vskip5mm
{\bf Corollary 2.1.} {\sl Every H\"offding kernel is positive definite:
For any collection $a_i, x_i \in \R$,
\be
\sum_{i,j=1}^n a_i a_j H(x_i,x_j) \geq 0.
\en
}

\vskip2mm
Usually, the equivalence of (2.2) and (2.3) is stated under the assumption 
that a kernel is continuous. In the case of H\"offding's kernels, 
this is however not important. Indeed, for the uniform distribution 
$U$ on the unit interval $(0,1)$, the corresponding kernel 
$H_U(x,y) = (x \wedge y)\,(1 - x \vee y)$ is positive definite 
on the square $0 \leq x,y \leq 1$. Hence the same is true for
$H_\mu(x,y) = H_U(F(x),F(y))$ on the plane.

Being positive definite, every H\"offding kernel satisfies
\be
H(x,y)^2 \leq H(x,x) H(y,y), \quad x,y \in \R,
\en
which may be used to construct a pseudometric
$$
d(x,y) = \big(H(x,x) - 2H(x,y) + H(y,y)\big)^{1/2}.
$$
This property can be strengthened in terms of the H\"offding measure
$\lambda$. Since, by the Cauchy inequality, 
$\cov_\mu(u,v)^2 \leq \Var_\mu(u) \Var_\mu(v)$, we get
from Theorem 1.1 that
$$
\bigg(\int_{-\infty}^\infty \int_{-\infty}^\infty 
f(x) g(y)\,d\lambda(x,y)\bigg)^2 \leq 
\int_{-\infty}^\infty \int_{-\infty}^\infty 
f(x) f(y)\,d\lambda(x,y) \ 
\int_{-\infty}^\infty \int_{-\infty}^\infty 
g(x) g(y)\,d\lambda(x,y)
$$
for all non-negative measurable
functions $f$ and $g$ on the real line. In particular, 
\be
\lambda(A \times B)^2 \, \leq \, 
\lambda(A \times A) \, \lambda(B \times B)
\en
for all Borel sets $A,B \subset \R$. Hence (2.4) appears
as an infinitesimal version of (2.5).

\vskip5mm
\section{{\bf Spectral decompositions}}
\setcounter{equation}{0}

\vskip2mm
\noindent
Since $H = H_\mu$ is positive definite, one may follow the advanced 
Mercer's theory on metric spaces and develop a canonical representation
\be
H(x,y) = \sum_{n=1}^\infty \alpha_n f_n(x) f_n(y)
\en
in terms of the eigenfunctions and eigenvalues of the linear operator
$$
T f(x) = \int_{-\infty}^\infty H(x,y)f(y)\,dy.
$$
Let
\be
a_0 = \inf\{x \in \R: F(x) > 0\}, \quad
a_1 = \sup\{x \in \R: F(x) < 1\},
\en
where $F$ is the distribution function of $X$. 
Applying Theorem 2.4 from Ferreira and Menegatto \cite{F-M} in the
setting of H\"offding's kernels, we obtain:

\vskip5mm
{\bf Corollary 3.1.} {\sl Suppose that $\mu$ has finite first absolute
moment, and $F$ is continuous on the interval $(a_0,a_1)$. There exists an orthonormal 
system of continuous functions $f_n$ in $L^2(a_0,a_1)$ and a non-increasing 
sequence $\alpha_n \geq 0$ such that $(3.1)$ holds for all $x,y \in (a_0,a_1)$.
This series is absolutely and uniformly convergent on finite proper subintervals
of $(a_0,a_1)$.
}

\vskip5mm
The moment assumption on $\mu$ guarantees that the H\"offding's kernel 
is square integrable over the rectangle $(a_0,a_1) \times (a_0,a_1)$. In particular,
$T$ is acting on the Hilbert space $L^2(a_0,a_1)$ as a compact and self-adjoint
operator. Moreover, it is a trace class so that
$$
{\rm tr}(T) = \int_a^b H(x,x)\,dx = \int_{-\infty}^\infty F(x)(1-F(x))\,dx =
\sum_{n=1}^\infty \alpha_n,
$$
where the series is convergent. The latter integral may also be recognized
as $\frac{1}{2}\,\E\,|X - X'|$ with $X'$ being an independent copy of $X$.

As a consequence of (2.6), the variance representation (2.1) may be expressed 
in the form
$$
\Var(u(X)) = \sum_{n=1}^\infty \alpha_n
\Big(\int_{a_0}^{a_1} u'(x) f_n(x) \,dx\Big)^2.
$$

For example, for the uniform distribution $\mu = U$, 
(2.6) holds for all $x,y \in (0,1)$ with $\alpha_n = 1/(n\pi)^2$ and 
$f_n(x) = \sqrt{2}\,\sin(n\pi x)$.

In the general situation, let $\mu$ have a continuous positive
density $p(x)$ in $a_0<x<a_1$. As easy to see, the spectral equation
$Tf = \alpha f$ is reduced to the Sturm-Liouville equation
$$
\alpha\,\Big(\frac{f'}{p}\Big)' + f = 0.
$$
When $p$ is continuous and positive on the finite interval
$[a_0,a_1]$, we thus arrive at the regular Sturm-Liouville problem with
boundary conditions $f(a_0) = f(a_1) = 0$, for which the spectral theory is
well-developed as well.

\vskip5mm
\section{{\bf Marginals of H\"offding measures}}
\setcounter{equation}{0}

\vskip2mm
\noindent
Since the kernel $H = H_\mu$ is symmetric about the diagonal $x=y$, the
H\"offding measure $\lambda = \lambda_\mu$ has equal marginals 
$\Lambda = \Lambda_\mu$ defined by
\be
\Lambda(A) = \lambda(A \times \R) = \int_A^\infty \int_{-\infty}^\infty
H(x,y)\,dx\,dy, \quad A \subset \R \ ({\rm Borel}).
\en
Obviously, it is absolutely continuous with respect to the Lebesgue measure
and is supported on the interval $(a_0,a_1)$, finite or not, defined in (3.2).
Concerning its density, let us emphasize the following two simple
properties.

\vskip5mm
{\bf Proposition 4.1.} {\sl If $X$ has finite first absolute moment, then 
the marginal $\Lambda$ is finite and has density
\be
h(x) = \frac{d\Lambda(x)}{dx} = \int_x^\infty (y-a)\,dF(y), \quad
a = \E X.
\en
In particular, it is unimodal with mode at the point $a$, that is, $h(x)$ is 
non-decreasing on the half-axis $x < a$ and is non-increasing for $x > a$.
Moreover, it is contintuous at $x=a$ with
\be
h(a-) = h(a+) = \frac{1}{2}\,\E\,|X-a|.
\en
If $\E\,|X| = \infty$, then the density of $\Lambda$
is a.e. infinite on $(a_0,a_1)$.
}

\vskip5mm
{\bf Proof.} According to (4.1), the measure $\Lambda$ has density
\begin{eqnarray}
h(x)
 & = &
\int_{-\infty}^\infty F(x \wedge y)\,(1 - F(x \vee y))\,dy \nonumber \\ 
 & = &
(1 - F(x)) \int_{-\infty}^x F(y)\,dy + F(x) \int_x^\infty (1-F(y))\,dy.
\end{eqnarray}
If $\E\,|X| = \infty$, then at least one of the two last integrals must be infinite
for all $a_0 < x < a_1$, which means that $h(x) = \infty$ a.e. on $(a_0,a_1)$.

If $\E\,|X| < \infty$, both integrals in (4.4) are finite. Assuming without loss
of generality that $F$ is continuous at the point $x$, one may integrate by parts
in (4.4) to obtain (4.2). Finally, since the function $x \rightarrow x-a$ is
vanishing at the point $a$, it follows that
$$
h(a+) = \int_a^\infty (y-a)\,dF(y) = \E\,(X-a)^+ = \frac{1}{2}\,\E\,|X-a|.
$$
A similar equality is also true for $h(a-)$, thus implying (4.3).
\qed

\vskip5mm
{\bf Proposition 4.2.} {\sl Assuming that $X$ has finite first absolute moment, 
the marginal $\Lambda$ is a multiple of $\mu$, if and only if $\mu$ is Gaussian.
}

\vskip5mm
{\bf Proof.} Without loss of generality we may assume that $\E X = 0$.
Since $\Lambda(\R) = \lambda(\R \times \R)$, the property $\Lambda = \sigma^2 \mu$ 
with some constant $\sigma^2$ implies that $\Lambda$ and $\lambda$ are finite and
forces $\mu$ to have a finite second moment. In that case, 
it follows from (1.6) that the Fourier-Stieltjes transform of 
$\Lambda$ is given by
$$
\hat \Lambda(t) = \hat \lambda(t,0) = 
\int_{-\infty}^\infty e^{itx}\,h(x)\,dx = -\frac{f'(t)}{t}, \quad 
t \in \R, \ t \neq 0,
$$
where $f$ is the characteristic function of $X$.
Hence, $\Lambda = \sigma^2 \mu$ if and only if 
$$
f'(t) = -\sigma^2 t f(t) \quad {\rm for \ all} \ \ t \in \R. 
$$
But this is only possible when $\mu$ is the Gaussian measure
with mean zero and variance $\sigma^2$.
\qed

\vskip5mm
Let us conclude with a few remarks. Often, the marginals of H\"offding's measures 
appear  in the particular case of the covariance representation (1.3) with 
the function $v(x) = x$. Then we arrive at
$$
\cov(X,u(X)) = \int_{-\infty}^\infty u'(x)\,h(x)\,dx,
$$
holding true as long as the integral is convergent.
If $\mu$ is supported on an interval $\Delta$ and has there an a.e.
positive density $p$, this formula may be rewritten as
\be
\cov(X,u(X)) = \E\,\tau(X) u'(X).
\en
Here the function
$$
\tau(x) = \frac{h(x)}{p(x)} =
\frac{1}{p(x)} \int_x^\infty (y-a)\,p(y)\,dy, \quad x \in \Delta,
$$
is called the Stein kernel. We have $\tau(x) = 1$ (a.e. on $\Delta$),
if and only if $\mu$ is the standard Gaussian measure $\gamma$, in which case
(4.5) becomes Stein's equation $\E\, X u(X) = \E\,u'(X)$. After the pioneering 
work \cite{St}, the identity (4.5) served as a starting point in the extensive development 
of Stein's method as an approach to various forms of the central limit 
theorem and estimating the distances to the normal law $\gamma$
avoiding the method of characteristic functions. For example,
Cacoullos, Papathansiou and Utev \cite{C-P-U} proposed a general upper bound
for the total variation distance
$$
\|\mu - \gamma\|_{\rm TV} \leq 4\,\E\,|\tau(X) - 1| + 4\,|1-a|
$$
and applied it in the proof of the CLT with respect to this strong distance.
For a comprehensive exposition of the whole theory, we refer an interested 
reader to the book by Chen, Goldstein and Shao \cite{C-G-S} and survey \cite{Sh}.

\vskip5mm
\section{{\bf Periodic covariance representations}}
\setcounter{equation}{0}

\vskip2mm
\noindent
A natural multidimensional extension of the covariance representation (1.3)
for a given probability measure $\mu$ on $\R^n$ could be the identity
\be
{\rm cov}_\mu(u,v) = \int_{\R^n} \int_{\R^n}
\left<\nabla u(x), \nabla v(y)\right> d\lambda(x,y),
\en
where $\lambda$ is a suitable measure on $\R^n \times \R^n$.
This is indeed possible when $\mu$ is Gaussian with covariance
matrix $\sigma^2 I$ (a multiple of the identity matrix). In this case,
$\lambda$ is unique and can be described in several equivalent 
ways including Ornstein-Uhlenbeck semigroups and interpolation
(\cite{Le}, \cite{H-P}). In \cite{B-G-H}, it was shown that the existence
of $\lambda$ in (5.1) forces the measure $\mu$ be Gaussian which gives
another characterization of this class in terms of covariance representations.
Nevertheless, some other variants of (5.1) could be applicable in order to involve
larger classes of probability distributions in such identities.
In particular, one can show that there is a spherical counterpart of (5.1),
\be
{\rm cov}_{\sigma_{n-1}}(u,v) = \int_{S^{n-1}} \int_{S^{n-1}}
\left<\nabla_S u(x), \nabla v_S(y)\right> d\lambda(x,y),
\en
with respect to the uniform distribution $\sigma_{n-1}$ on the unit sphere
$S^{n-1}$ in $\R^n$ for some specific measure $\lambda$ on $S^{n-1} \times S^{n-1}$ 
having multiples of $\sigma_{n-1}$ as marginals. Here
$\nabla_S u$ denotes the spherical gradient of a smooth function $u$ on the sphere.

However, the measures $\lambda$ in such representations are not unique
anymore. This can be seen already in the case of the circle $S^1$, when (5.2)
is reduced to the covariance representation
\be
{\rm cov}_\mu(u,v) = \int_0^1 \int_0^1 u'(x) v'(y)\, d\lambda(x,y)
\en
for the uniform distribution $\mu = m$ on $(0,1)$ in the class of all
1-periodic smooth functions $u$ and $v$ on the real line.
Here, $\lambda$ is a certain finite measure on the square which we allow to be
a signed measure for the sake of generality. Without loss of generality,
we require that it is symmetric about the diagonal $x=y$.

Keeping aside the multidimensional setting for a separate consideration,
in what follows we focus on (5.3), assuming that $\mu$ is a given probability 
measure on $[0,1)$.

\vskip5mm
{\bf Definition.} We call a signed symmetric Borel measure $\lambda$ 
on $[0,1) \times [0,1)$ a mixing measure for $\mu$, if (5.3) holds true for all
1-periodic smooth functions $u$ and $v$ on the real line.

\vskip5mm
As we discussed before, this identity always holds with the H\"offding measure
$\lambda = \lambda_\mu$. But, its marginals may be a multiple of $\mu$
in the Gaussian case only. This motivates the following:

\vskip5mm
{\bf Question.} Given $\mu$, is it possible to choose a mixing measure $\lambda$
whose marginals are multiples of $\mu$? If so, how to describe all of them
and choose a best one (in some sense)?

\vskip5mm
Towards this question, we prove the next assertion.

\vskip5mm
{\bf Theorem 5.1.} {\sl Let $\mu$ be a probability measure on $[0,1)$ with 
the H\"offding measure $\lambda_\mu$. Subject to the constraint that the 
marginal distribution of $\lambda$ in $(5.3)$ is equal to $c\mu$ for 
a prescribed value $c \in \R$, the mixing measure $\lambda$ exists, is unique, 
and is given by
\be
\lambda = \lambda_\mu + (\sigma^2 - c)\,m \otimes m + c\, 
(\mu \otimes m + m \otimes \mu) -
(\Lambda_\mu \otimes m + m \otimes \Lambda_\mu),
\en
where $\Lambda_\mu$ is the marginal of $\lambda_\mu$ 
and $\sigma^2$ is the variance of $\mu$.
}

\vskip5mm
\section{{\bf Covariance representations for the uniform distribution}}
\setcounter{equation}{0}

\vskip2mm
\noindent
We postpone the proof of Theorem 5.1 to Sections 8-9. The most interesting 
case in the periodic representation (5.3) is the one where $\mu$ represents
a uniform distribution $m$ on (0,1). Let us specialize Theorem 5.1 
to this case and consider identities of the form
\be
{\rm cov}_m(u,v) = \int_0^1 \int_0^1 u'(x) v'(y)\, d\lambda(x,y).
\en
As a consequence of Theorem 5.1, we obtain the following statement needed
in the study of covariance identities on the circle (as part of multidimensional
spherical identities (5.2)).

\vskip5mm
{\bf Corollary 6.1.} {\sl Subject to the constraint that the marginal 
distribution of a mixing measure $\lambda$ in $(6.1)$ is equal 
to $cm$, $c \in \R$, the measure $\lambda$ is unique and has density
\be
\frac{\lambda(x,y)}{dx\,dy} = D(|x-y|) + \Big(c - \frac{1}{24}\Big), 
\quad x,y \in (0,1),
\en
where
\be
D(h) = \frac{1}{8}\,\big[\,1 - 4h(1 - h)\big], \quad 0 \leq h \leq 1.
\en
}

\vskip2mm
Note that $D(h) \geq 0$ for all $h \in [0,1]$ and the inequality becomes
an equality for $h=\frac{1}{2}$.

Moreover, the mixing measure $\lambda$ is non-negative, if and only if 
$c \geq \frac{1}{24}$. Hence, the smallest 
positive measure (for the usual comparison) corresponds to the parameter 
$c = \frac{1}{24}$, when it has the density $\psi(x,y) = D(|x-y|)$. In this sense, 
the optimal variant of (6.1) is given by the covariance representation
\bee
{\rm cov}_{m}(u,v)
 & = &
\int_0^1 \int_0^1 u'(x) v'(y)\,D(|x-y|)\,dx\,dy \\
 & = &
\frac{1}{24} \int_0^1 \int_0^1 u'(x) v'(y)\,d\nu(x,y)
\ene
with a probability measure $d\nu(x,y) = 24\,D(|x-y|)\,dx\,dy$
on $(0,1) \times (0,1)$.
It has the uniform distribution $m$ on $(0,1)$ as a marginal one.

\vskip5mm
{\bf Proof.} Returning to (5.4), note that, if $\mu$ has density $p$, 
then the measure $\lambda$ has density
\be
\psi(x,y) = h(x,y) + 
(\sigma^2 - c) + c\,(p(x) + p(y)) - (q(x) + q(y))
\en
on $[0,1) \times [0,1)$, where $h$ denotes the density of the marginal
$\Lambda_\mu$ of the H\"offding measure $\lambda_\mu$ with density
$H(x,y) = F(x \wedge y)\,(1 - F(x \vee y))$. Recall that according to (4.1), 
if the mean of $\mu$ is $a = \E X$, then
$$
h(x) = \int_x^1 y p(y)\,dy - a(1-F(x)), \quad 0 \leq x \leq 1.
$$

In the case of the uniform distribution $\mu = m$, the distribution 
function and density are given by $F(x) = x$ and $p(x) = 1$, $0 \leq x \leq 1$. 
Then, by (6.4),
\be
\psi(x,y) = H(x,y) + (\sigma^2 + c) - (h(x) + h(y))
\en
with $H(x,y) = (x \wedge y)\,(1 - (x \vee y))$ and
$\sigma^2 = \frac{1}{12}$. To simplify, assume that 
$0 \leq x \leq y \leq 1$. Then
\bee
\psi(x,y) 
 & = &
\Big(c + \frac{1}{12}\Big) + x (1 - y) - \frac{1}{2}\,((x - x^2) + (y - y^2)) \\
 & = &
\Big(c + \frac{1}{12}\Big) -\frac{1}{2}\,\Big((y-x)(1 - (y-x))\Big) \\
 & = &
\Big(c - \frac{1}{24}\Big) + D(|x-y|).
\ene
\qed

\vskip2mm
{\bf Remark.} If we want to write down a similar representation on 
the interval $(0,T)$, $T>0$, one may use a linear transform. 
Let $m_T$ denote the uniform distribution on $(0,T)$. Then we get that, 
for all smooth $T$-periodic functions $u$ and $v$,
and for all $c \geq 1/24$, 
$$
{\rm cov}_{m_T}(u,v) = \int_0^T \int_0^T u'(x) v'(y)\,d\lambda_T(x,y)
$$
with a positive measure having the density
$$
\frac{d\lambda_T(x,y)}{dx\,dy} = 
D\Big(\Big|\frac{x}{T} - \frac{y}{T}\Big|\Big) + \Big(c - \frac{1}{24}\Big)
$$
on $(0,T) \times (0,T)$. It has the marginal $cT m_T(dx) = c\,dx$ on $(0,T)$.

\vskip5mm
\section{{\bf Densities bounded away from zero}}
\setcounter{equation}{0}

\vskip2mm
\noindent
In the general situation, the question of whether or not the mixing measure
$\lambda$ is positive for a certain value of $c$ is rather interesting 
(in which case this constant has to be positive as well). Here we give one 
sufficient condition generalizing the previous example of the uniform distribution. 
As usual, we denote by $\sigma$ is the standard deviation of a random variable 
$X$ distributed according to $\mu$.

\vskip5mm
{\bf Corollary 7.1.} {\sl Suppose that the probability measure $\mu$ on $(0,1)$
has a density $p$ such that $p(x) \geq \alpha$ for all $x \in (0,1)$ with
some constant $\alpha > \frac{1}{2}$. Then there exists a positive mixing
measure $\lambda$ in the periodic covariance representation
\be
{\rm cov}_\mu(u,v) = \int_0^1 \int_0^1 u'(x) v'(y)\, d\lambda(x,y),
\en
whose marginal is a multiple $c\mu$ of $\mu$. One may choose
$c = \frac{\sigma - \sigma^2}{2\alpha - 1}$.
}

\vskip5mm
{\bf Proof.} According to (5.4), subject to the constraint that the marginal 
distribution of $\lambda$ in $(7.1)$ is equal to $c\mu$, the mixing measure $\lambda$
has density
$$
p(x,y) = H(x,y) + \sigma^2 + c\,(p(x) + p(y) - 1) - (h(x) + h(y)), \quad
x,y \in (0,1),
$$
where $H(x,y)$ is the H\"offding kernel and $h(x)$ is the density of the mariginal
distribution $\Lambda$. Hence, it is non-negative, as long as
$$
c\,(p(x) + p(y) - 1) \geq h(x) + h(y) - \sigma^2.
$$
By the assumption, $p(x) + p(y) - 1 \geq 2\alpha -1$, so that,
it is sufficient to require that
\be
c\,(2\alpha - 1) \geq h(x) + h(y) - \sigma^2.
\en

Now, let us recall that, by Proposition (4.1), $h(x)$ is unimodal and
continuous. Morever, according to (4.3), for all $x \in (0,1)$,
$$
2h(x) \leq 2h(a) = \E\,|X-a| \leq \sigma, \quad a = \E X,
$$
where we applied Cauchy's inequality. Note that $\sigma^2 < \sigma$.
Hence, the right-hand side of (7.2) is bounded from above by $\sigma - \sigma^2$.
\qed

\vskip5mm
{\bf Example.} The symmetric beta distribution with parameters $(\frac{1}{2},\frac{1}{2})$,
that is, with density
$$
p(x) = \frac{1}{\pi \sqrt{x(1-x)}}, \quad 0 < x < 1,
$$
satisfies the conditions of Corollary 7.1 with $\alpha = \frac{2}{\pi}$ and
$\sigma^2 = \frac{1}{8}$. Hence, the conclusion in this corollary is true with
$c = \frac{\pi}{8 (4 - \pi)}\,(\sqrt{8} - 1) \sim 0.8364...$. In fact, a more careful
analysis shows that one may take $c = \frac{1}{8}$ in this example
(which is optimal).

\vskip5mm
\section{{\bf Characterization of mixing measures}}
\setcounter{equation}{0}

\vskip2mm
\noindent
Let us first comment on the uniqueness issue in the problem
of characterization of mixing measures in the covariance representation
(5.3). Applying this identity to the exponential functions 
$u(x) = e^{2\pi ikx}$ and $v(y) = e^{2\pi ily}$, we get the relation
\be
f(k+l) - f(k) f(l) = -(2\pi)^2\,kl\,\hat \lambda(k,l)
\en
for all integers $k,l$, where
$$
\hat \lambda(k,l) = \int_0^1 \int_0^1
e^{2\pi i(kx + ly)}\,d\lambda(x,y), \quad k,l \in \Z,
$$
denotes the Fourier transform of $\lambda$ restricted to integers. 
By the Stone-Weierstrass theorem (applied on the circle), $\hat \lambda$ 
determines any signed Borel measure $\lambda$ on $[0,1) \times [0,1)$ 
in a unique way. This transform is explicitly defined in (8.1) as long as 
$k,l \neq 0$. Otherwise, both sides of (8.1) are vanishing. 
The fact that (8.1) does not define $\hat \lambda$ for all integers 
does not allow us to reconstruct $\lambda$.

Moreover, due to the periodicity of $u$ and $v$, we have
$$
\int_0^1 \int_0^1 u'(x) v'(y)\,d\Lambda_1(x)\,dm(y) =
\int_0^1 \int_0^1 u'(x) v'(y)\,dm(x)\,d\Lambda_2(y) = 0
$$
for all signed measures $\lambda_1$ and $\lambda_2$ on $[0,1)$, where 
$m$ denotes the uniform probability measure on that interval. Hence, 
once (5.3) is fulfilled for a measure $\lambda$, in particular, for 
the H\"offding measure $\lambda_\mu$, it is also fulfilled for
\be
\lambda = \lambda_\mu + \Lambda_1 \otimes m + m \otimes \Lambda_2
\en
for any choice of signed measures $\Lambda_1$ and $\Lambda_2$ on $[0,1)$.
We also have the converse statement (where the symmetry requirement is
not required for a moment).

\vskip5mm
{\bf Lemma 8.1.} {\sl Let $\mu$ be a Borel probability measure on $[0,1)$.
The covariance representation 
\be
{\rm cov}_\mu(u,v) = \int_0^1 \int_0^1 u'(x) v'(y)\, d\lambda(x,y)
\en
holds true
for all smooth, 1-periodic functions $u$ and $v$,
if and only if $\lambda$ has the form $(8.2)$ for some (arbitrary)
signed measures $\Lambda_1$ and $\Lambda_2$ on $[0,1)$.
}

\vskip5mm
{\bf Proof.}  We only need to consider the necessity part.
Assume that (8.3) holds true for all $C^1$-smooth periodic functions 
$u$ and $v$, so that
\be
\int_0^1 \int_0^1 u'(x) v'(y)\,d\lambda_\mu(x,y) =
\int_0^1 \int_0^1 u'(x) v'(y)\,d\lambda(x,y).
\en
Putting $f = u'$, $g = v'$, we then have
\be
\int_0^1 \int_0^1 f(x) g(y)\,d\lambda_\mu(x,y) =
\int_0^1 \int_0^1 f(x) g(y)\,d\lambda(x,y).
\en
Due to the periodicity of $u$ and $v$, necessarily
\be
\int_0^1 f(x) \,dm(x) = \int_0^1 g(y)\,dm(y) = 0
\en
and
\be
f(0) = f(1), \quad  g(0) = g(1).
\en
Conversely, starting from continuous $f$ an $g$ 
on $[0,1]$ satisfying (8.6)-(8.7), we may define the functions
$$
u(x) = \int_0^x f(t)\,dt, \quad
v(y) = \int_0^y g(s)\,ds,
$$
which have $C^1$-smooth 1-periodic extensions 
from $[0,1)$ to the whole real line and satisfy (8.4).
Thus, our hypothesis (8.3) is equivalent to (8.5) subject to 
(8.6)--(8.7).

Let us reformulate the latter by identifying $[0,1)$ with the circle 
$S^1$ via the map $x \rightarrow e^{2\pi ix}$. It pushes forward $m$ 
to the uniform probability measure $\sigma_1$ on the circle and pushes 
$\lambda - \lambda_\mu$ to some signed
measure $L$ on the torus $S^1 \times S^1$. Hence, (8.5) 
subject to (8.6)--(8.7) is the same as the requirement
\be
\int_{S^1} \int_{S^1} \xi(t) \eta(s)\,dL(t,s) = 0
\en
in the class of all continuous functions $\xi,\eta$ on the circle
such that $\int_{S^1} \xi\,d\sigma_1 = \int_{S^1} \eta\,d\sigma_1 = 0$.
The latter assumption may be dropped, if we rewrite (8.8) as
\be
\int_{S^1} \int_{S^1} (\xi(t) - \bar \xi)\, (\eta(s) - \bar \eta)\,dL(t,s) = 0,
\en
where
$$
\bar \xi = \int_{S_1} \xi\,d\sigma_1, \quad 
\bar \eta = \int_{S_1} \eta\,d\sigma_1.
$$
In this step, (8.9) is readily extended to the class of all bounded,
Borel measurable functions $\xi$ and $\eta$ on $S^1$.
Using the marginal measures
$$
L_1(A) = L(A \times S^1), \quad L_2(B) = L(S^1 \times B),
$$
one may now rewrite the equality (8.9) as
\bee
\int \!\!\int \xi(t) \eta(s)\,dL(t,s) 
 & = & 
\int \xi(t)\,d\sigma_1(t) \int \eta(s) \,d\sigma_1(s) \\
 & & \hskip-10mm + \
\int \xi(t)\,d\sigma_1(t) \int \eta(s)\,dL_2(s) + 
\int \xi(t)\,dL_1(t) \int \eta(s)\,d\sigma_1(s).
\ene
But this means that
$$
L = \sigma_1 \otimes \sigma_1 + \sigma_1 \otimes L_2 + 
L_1 \otimes \sigma_1.
$$
Pushing $\sigma_1$ and $L_j$ back to $[0,1)$ with images $\Lambda_j$
and $L$ to $[0,1) \times [0,1)$, we arrive at 
$$
\lambda - \lambda_\mu = m \otimes m + m \otimes \Lambda_2 + 
\Lambda_1 \otimes m.
$$
This is an equivalent form for (8.2).
\qed

\vskip5mm
\section{{\bf Proof of Theorem 5.1}}
\setcounter{equation}{0}

\vskip2mm
\noindent
Recall that the mixing measure $\lambda$ in (5.3) has to be supported 
on $[0,1) \times [0,1)$ and is required to be symmetric about the diagonal 
of this square. This is fulfilled for the H\"offding measure $\lambda_\mu$. 
Hence, the measure $\lambda$ of the form (8.2) is symmetric about the diagonal 
of the square, if and only if $\Lambda_2 - \Lambda_1$ is proportional 
to the uniform measure $m$. 
In other words, the class of all symmetric measures $\lambda$ satisfying 
the covariance representation (5.3) is described by the formula
\be
\lambda = \lambda_\mu + b\,m \otimes m + \Lambda \otimes m + 
m \otimes \Lambda
\en
with arbitrary $b \in \R$ and arbitrary signed measures $\Lambda$ on $[0,1)$. 

Such measures have equal marginals
$$
{\rm Proj}(\lambda) = \Lambda_\mu + (b+Q)\,m + \Lambda, \quad
\Lambda_\mu = {\rm Proj}(\lambda_\mu),
$$
where $Q = \Lambda([0,1))$ and $\Lambda_\mu$ is the marginal of 
$\lambda_\mu$ described in (4.1). We want this measure 
to be a multiple of the original probability measure $\mu$ on $[0,1)$, that is,
$$
c\mu = \Lambda_\mu + (b+Q)\,m + \Lambda
$$
for some prescribed value $c \in \R$. Then, necessarily with some $d \in \R$
$$
\Lambda = c\mu - \Lambda_\mu + d m
$$
To determine the value of $d$, we plug this into (9.1) and get
\begin{eqnarray}
\lambda 
 & = &
\lambda_\mu + b\,m \otimes m + (c\mu - \Lambda_0 + d m) \otimes m + 
m \otimes (c\mu - \Lambda_\mu + d m) \nonumber \\
 & = &
\lambda_\mu + (b+2d)\,m \otimes m + c\, 
(\mu \otimes m + m \otimes \mu) -
(\Lambda_\mu \otimes m + m \otimes \Lambda_\mu).
\end{eqnarray}
On marginals this equality becomes the relation 
\be
b + c + 2d - \sigma^2 = 0,
\en
where 
$$
\sigma^2 = \Lambda_\mu([0,1)) = \lambda_\mu([0,1) \times [0,1))
$$
is variance of a random variable distributed according to $\mu$.
This is how $d$ should be determined in terms of the
free parameters $b$ and $c$. In that case,
\be
\Lambda = c\mu - \Lambda_\mu + \frac{\sigma^2 - b - c}{2}\, m.
\en
It remains to apply (9.3) in (9.2). Theorem 5.1 is now proved.

\vskip5mm
\section{{\bf Proof of Theorem 1.1}}
\setcounter{equation}{0}

\vskip2mm
\noindent
One may assume that the locally absolutely functions $u$ and $v$ 
are real-valued and have Borel measurable Radon-Nikodym derivatives
$u'$ and $v'$ (which are locally integrable). 

{\bf Step 1.} Assume that $u$ and $v$ have
non-negative $u'$ and $v'$. In view of the monotonicity of 
$u$ and $v$, the covariance of $u(X)$ and $v(X)$ is well-defined
and is given by
\bee
\cov(u(X),v(X)) 
 & = &
\int\!\! \int_{x < y} 
(u(x) - u(y))(v(x) - v(y))\,d\mu(x)\,d\mu(y) \\
 & = &
\int_{-\infty}^\infty \int_{-\infty}^\infty
\int_{-\infty}^\infty \int_{-\infty}^\infty
u'(t) v'(s)\, 1_{\{x \leq t < y\}}\, 1_{\{x \leq s < y\}}
\,dt\,ds\,d\mu(x)\,d\mu(y)
\ene
as the Lebesgue integral on $\R^4$ over the product measure 
$L \otimes L \otimes \mu \otimes \mu$
(where $L$ is the Lebesgue measure on the real line).
By Fubini's theorem, we obtain that
\be
{\rm cov}(u(X),v(X)) = 
\int_{-\infty}^\infty \int_{-\infty}^\infty 
u'(x) v'(y)\,H(x,y)\,dx\,dy.
\en
In particular,
\be
\Var(u(X)) = \int_{-\infty}^\infty \int_{-\infty}^\infty 
u'(x) u'(y)\,H(x,y)\,dx\,dy,
\en
and similarly for $v$. 
As a by-product, using
$|{\rm cov}(u(X),v(X))|^2 \leq \Var(u(X)) \, \Var(v(X))$,
\begin{eqnarray}
\bigg(\int_{-\infty}^\infty \int_{-\infty}^\infty 
u'(x) v'(y)\,H(x,y)\,dx\,dy\bigg)^2
 \leq & & \nonumber \\
 & & \hskip-65mm
\int_{-\infty}^\infty \int_{-\infty}^\infty 
u'(x) u'(y)\,H(x,y)\,dx\,dy \ 
\int_{-\infty}^\infty \int_{-\infty}^\infty 
v'(x) v'(y)\,H(x,y)\,dx\,dy.
\end{eqnarray}

{\bf Step 2.} 
In the general case, define locally absolutely continuous,
non-decreasing functions
$$
\tilde u(x) = \int_0^x |u'(t)|\,dt, \quad
\tilde v(x) = \int_0^x |v'(t)|\,dt,
$$
which have Radon-Nikodym derivatives $|u'|$ and $|v'|$. Since
$|\tilde u(x) - \tilde u(y)| \geq |u(x) - u(x)|$
for all $x,y \in \R$, we have
$$
\Var(\tilde u(X)) \geq \Var(u(X)),
$$
and similarly for $v$. By the previous step, 
$$
\Var(\tilde u(X)) = 
\int_{-\infty}^\infty \int_{-\infty}^\infty 
|u'(x)|\, |u'(y)|\,H(x,y)\,dx\,dy,
$$
and the same is true for $\tilde v$. Since this and a similar integral 
for $v$ are supposed to be finite, we conclude that both $u(X)$ and $v(X)$ 
have finite second moments.

One may now repeat the arguments from Step 1, using the inequality (10.3) 
with $|u'|$ and $|v'|$ in place of $u'$ and $v'$ respectively. This will 
justify an application of the Fubini's theorem, and then we obtain 
the identity (10.1) under the conditions in (1.5). This also insures
the integrability of $u'(x) v'(y)$ over $\lambda$
as a consequence of (1.5) and (10.2).

{\bf Step 3.} For the uniqueness issue, let $\lambda$ be a locally finite 
measure on the plane satisfying (1.3) in the class
$C_b^\infty$. Using a simple approximation, we obtain that
$\lambda(A \times B) = \lambda_\mu(A \times B)$ for all bounded intervals
$A$ and $B$. Hence, this equality is true for all Borel subsets of $\R^2$.

\vskip5mm
{\bf Remarks.} In \cite{Lo}, Theorem 1.1 is proved in a more general
setting of random variables $X$ and $Y$ as the identity (1.1),
assuming that $u$ and $v$ are absolutely continuous (not just locally), that is,
\be
\int_{-\infty}^\infty |u'(x)|\,dx < \infty, \quad 
\int_{-\infty}^\infty |v'(x)|\,dx < \infty,
\en
and such that $u(X)$, $v(Y)$, $u(X) v(Y)$ have finite first absolute moments
(cf. Theorem 3.1 in \cite{Lo}). Note, however, that the condition (10.4) insures 
that both $u$ and $v$ are bounded, so that the moment assumptions are 
fulfilled automatically. A similar assertion with $X=Y$ is given in \cite{S-W1},
Corollary 4, where in addition to the absolute continuity it is assumed that
\be
\E\,|u(X)|^p < \infty, \quad \E\,|v(X)|^q < \infty
\en
for some $p,q \geq 1$ such that $\frac{1}{p} + \frac{1}{q} = 1$. Again,
the latter assumption is not needed, if we assume (10.4). As for
the more general case of locally absolutely continuous $u$ and $v$,
the condition (10.5) and even the assumption on the boundedness of these
functions do not guarantee that the integral in (1.3) is convergent in the
Lebesgue sense, that is,
$$
\int_{-\infty}^\infty \int_{-\infty}^\infty 
|u'(x)|\, |v'(y)|\,H_\mu(x,y)\,dx\,dy < \infty.
$$
For example, for $u(x) = v(x) = \cos x$, this integral is divergent as long as
$\E\, |X| = \infty$.


\end{document}